\documentclass[11pt]{amsart}
\usepackage{amsmath,amsfonts}

\newtheorem{thm}{Theorem}[section]
\newtheorem{cor}[thm]{Corollary}
\newtheorem{lem}[thm]{Lemma}
\newtheorem{fact}[thm]{Fact}
\newtheorem{prop}[thm]{Proposition}
\newtheorem{rem}[thm]{Remark}

\newtheorem{notation}[thm]{Notation}
\numberwithin{equation}{section}

\begin{document}

\title[Smooth approximation by mappings with no critical points]
{Uniform approximation of continuous mappings by smooth
mappings with no critical points on Hilbert manifolds}
\author{Daniel Azagra and Manuel Cepedello Boiso}
\date{March 22, 2002}
\thanks{The first-named author was supported by a
Marie Curie Fellowship of the European Community Training and
Mobility of Researchers Programme under contract number
HPMF-CT-2001-01175}

\begin{abstract}
We prove that every continuous mapping from a separable
infinite-dimensional Hilbert space $X$ into $\mathbb{R}^{m}$ can
be uniformly approximated by $C^\infty$ smooth mappings {\em with
no critical points}. This kind of result can be regarded as a
sort of very strong approximate version of the Morse-Sard
theorem. Some consequences of the main theorem are as follows.
Every two disjoint closed subsets of $X$ can be separated by a
one-codimensional smooth manifold which is a level set of a
smooth function with no critical points; this fact may be viewed
as a nonlinear analogue of the geometrical version of the
Hahn-Banach theorem. In particular, every closed set in $X$ can
be uniformly approximated by open sets whose boundaries are
$C^\infty$ smooth one-codimensional submanifolds of $X$. Finally,
since every Hilbert manifold is diffeomorphic to an open subset
of the Hilbert space, all of these results still hold if one
replaces the Hilbert space $X$ with any smooth manifold $M$
modelled on $X$.
\end{abstract}

\maketitle

\section[Introduction]{Introduction and main results}

A fundamental result in differential topology and analysis is the
Morse-Sard theorem \cite{Sard1, Sard2}, which states that if
$f:\mathbb{R}^{n}\longrightarrow \mathbb{R}^{m}$ is a $C^r$ smooth
function, with $r>\max\{n-m, 0\}$, and $C_{f}$ stands for the set
of critical points of $f$ (that is, the points $x$ at which the
differential $df(x)$ is not surjective), then the set of critical
values, $f(C_{f})$, is of (Lebesgue) measure zero in
$\mathbb{R}^{m}$. This result also holds true for smooth functions
$f:X\longrightarrow Y$ between two smooth manifolds of dimensions
$n$ and $m$ respectively.

Several authors have dealt with the question as to what extent
one can obtain a similar result for infinite-dimensional spaces
or manifolds modelled on such spaces. Let us recall some of their
results.

Smale \cite{Smale} proved that if $X$ and $Y$ are separable
connected smooth manifolds modelled on Banach spaces and
$f:X\longrightarrow Y$ is a $C^r$ Fredholm map (that is, every
differential $df(x)$ is a Fredholm operator between the
corresponding tangent spaces) then $f(C_{f})$ is meager, and in
particular $f(C_{f})$ has no interior points, provided that
$r>\max\{\textrm{index}(df(x)), 0\}$ for all $x\in X$; here
index($df(x)$) stands for the index of the Fredholm operator
$df(x)$, that is, the difference between the dimension of the
kernel of $df(x)$ and the codimension of the image of $df(x)$,
which are both finite. However, these assumptions are quite
restrictive: for instance, if $X$ is infinite-dimensional then
there is no Fredholm map $f:X\longrightarrow\mathbb{R}$. In
general, the existence of a Fredholm map $f$ from a manifold $X$
into another manifold $Y$ implies that $Y$ is
infinite-dimensional whenever $X$ is.

On the other hand, one cannot dream of extending the Morse-Sard
theorem to infinite dimensions without imposing strong
restrictions. Indeed, as shown by Kupka's counterexample
\cite{Kupka}, there are $C^\infty$ smooth functions
$f:X\longrightarrow\mathbb{R}$, where $X$ is a Hilbert space, so
that their sets of critical values $f(C_{f})$ contain intervals
and in particular have non-empty interior.

More recently, S. M. Bates has carried out a deep study
concerning the sharpness of the hypothesis of the Morse-Sard
theorem and the geometry of the sets of critical values of smooth
functions. In particular he has shown that the above $C^r$
smoothness hypothesis in the statement of the Morse-Sard theorem
can be weakened to $C^{r-1,1}$; see \cite{Bates1, Bates2, Bates3,
Bates4, Bates5}. C. G. Moreira and Bates have studied some
generalizations of the Morse-Sard theorem related to Hausdorff
measures and Hausdorff dimensions. They have also shown that the
function $f$ as in Kukpa's counterexample can even be assumed to
be a polymonial of degree three; see \cite{Bates-Moreira,
Moreira}.

Nevertheless, for many applications of the Morse-Sard theorem, it
is often enough to know that any given continuous function can be
uniformly approximated by a map whose set of critical values has
empty interior. In this direction, Eells and McAlpin established
the following theorem \cite{EellsMcAlpin}: if $X$ is a separable
Hilbert space, then every continuous function from $X$ into
$\mathbb{R}$ can be uniformly approximated by a smooth function
$f$ whose set of critical values $f(C_{f})$ is of measure zero.
This allowed them to deduce a version of this theorem for
mappings between smooth manifolds $M$ and $N$ modelled on $X$ and
a Banach space $F$ respectively, which they called an {\em
approximate Morse-Sard theorem}: every continuous mapping from
$M$ into $N$ can be uniformly approximated by a smooth function
$f:M\longrightarrow N$ so that $f(C_{f})$ has empty interior.
However, this seemingly much more general version of the result
is a bit tricky: indeed, as they already observed
(\cite{EellsMcAlpin}, Remark 3A), when $F$ is
infinite-dimensional, the function $f$ they obtain satisfies that
$C_{f}=M$, although $f(M)$ has empty interior in $N$.
Unfortunately, even though all the results of that paper seem to
be true, some of the proofs are not correct.

In this paper we will prove a much stronger result: if $M$ is a
$C^\infty$ smooth manifold modelled on a separable
infinite-dimensional Hilbert space $X$ (in the sequel such a
manifold will be called a Hilbert manifold), then every continuous
mapping from $M$ into $\mathbb{R}^{m}$ can be uniformly
approximated by $C^\infty$ smooth mappings {\em with no critical
points}. This kind of result might be regarded as the strongest
possible one of any class of approximate Morse-Sard theorems,
when the target space is finite-dimensional.

As a by-product we also obtain the following: for every open set
$U$ in a separable Hilbert manifold $M$ there is a $C^{\infty}$
smooth function $f$ whose support is the closure of $U$ and so
that $d f(x)\neq 0$ for every $x\in U$. This result could be
summed up by saying that for every open subset $U$ of $M$ there
is a function $f$ whose open support is $U$ and which does not
satisfy Rolle's theorem; one should compare this result with the
main theorem from \cite{AJ} (see also the references therein).

Either of these results has in turn interesting consequences
related to smooth approximation and separation of closed sets. For
instance, every two disjoint closed subsets in $M$ can be separated 
by a smooth one-codimensional submanifold of $M$ which is a level set 
of a smooth function with no critical points. This may be
regarded as a nonlinear analogue of the geometrical version of the
Hahn-Banach theorem. In particular, every closed subset of $M$ can be
uniformly approximated by open sets whose boundaries are smooth 
one-codimensional submanifolds of $M$.

So far these are some {\em good} consequences of our main result,
all of them somehow related to Morse-Sard type theorems. But
there are some {\em bad} consequences as well, perhaps the most
noticeable one being that, since the set of smooth functions with
no critical points is dense in the set of continuous functions
defined on a Hilbert manifold, there are quite large sets of
smooth functions for which no conceivable Morse theory could be
valid.

Let us now formally state our main results. For the sake of a
convenient notation in our proofs, when $\varphi$ takes real
values we indistinctly use the symbols $d\varphi(x)=\varphi'(x)$
to denote the derivative of $\varphi$ at a point $x$, and we
reserve $d\varphi(x)$ for the derivative of a vector-valued
function $\varphi:M\longrightarrow\mathbb{R}^{m}$ at a point
$x\in M$.

\begin{thm}\label{main theorem}
Let $U$ be an open subset of a separable infinite-dimensional
Hilbert space $X$. Then, for every continuous mapping
$f:U\longrightarrow\mathbb{R}^{m}$ and for every continuous
positive function $\varepsilon:U\longrightarrow (0,+\infty)$,
there exists a $C^\infty$ smooth mapping
$\psi:U\longrightarrow\mathbb{R}^{m}$ such that
$\|f(x)-\psi(x)\|\leq\varepsilon(x)$ and $d\psi(x)$ is surjective
for all $x\in X$ (that is, $\psi$ has no critical points).
\end{thm}

We will prove this result in the following section. Let us now
establish the announced consequences of Theorem \ref{main
theorem}.

One could adapt the ideas in our proof to extend Theorem
\ref{main theorem} to the setting of Hilbert manifolds but, for
simplicity, we will instead use another approach. Indeed, bearing
in mind a fundamental result on Hilbert manifolds due to Eells
and Elworthy \cite{EE} that every separable Hilbert manifold can
be $C^\infty$ embedded as an open subset of the Hilbert space, it
is a triviality to observe that Theorem \ref{main theorem} still
holds if we replace $U$ with a a separable Hilbert manifold.

\begin{thm}\label{main theorem for manifolds}
Let $M$ be a separable Hilbert manifold. Then, for every
continuous mapping $f:M\longrightarrow\mathbb{R}^{m}$ and every
continuous positive function $\varepsilon:M\longrightarrow
(0,+\infty)$, there exists a $C^\infty$ smooth mapping
$\psi:M\longrightarrow\mathbb{R}^{m}$ such that $\psi$ has no
critical points and $\|f(x)-\psi(x)\|\leq\varepsilon(x)$ for all
$x\in X$.
\end{thm}
\begin{proof}
According to the main theorem of \cite{EE}, there is a $C^\infty$
embedding of $M$ onto an open subset of the Hilbert space $X$.
Therefore $M$ is $C^\infty$ diffeomorphic to an open subset $U$
of $X$; let $h:U\longrightarrow M$ be such a $C^\infty$
diffeomorphism. Consider the continuous mappings $g=f\circ
h:U\longrightarrow\mathbb{R}^{m}$ and $\delta=\varepsilon\circ
h:U\longrightarrow (0,+\infty)$. By Theorem \ref{main theorem}
there is a $C^\infty$ smooth function
$\varphi:U\longrightarrow\mathbb{R}^{m}$ so that $\varphi$ has no
critical points, and
    $$
    \|g(y)-\varphi(y)\|\leq\delta(y)
    $$
for all $y\in U$. Now define $\psi=\varphi\circ
h^{-1}:M\longrightarrow\mathbb{R}^{m}$. Since $h$ is a
diffeomorphism it is clear that $h$ takes the critical set of
$\psi$ onto the critical set of $\varphi=\psi\circ h$. But, as
the latter is empty, so is the former; that is, $\psi$ has no
critical points either. On the other hand, we have that
    $$
    \|f(x)-\psi(x)\|=\|g(h^{-1}(x))-\varphi(h^{-1}(x))\|\leq
    \delta(h^{-1}(x))=\varepsilon(x)
    $$
for all $x\in M$.
\end{proof}

As an easy corollary we can deduce our promised nonlinear version
of the geometrical Hahn-Banach theorem.

We will say that an open subset $U$ of a Hilbert manifold $M$ is
{\em smooth} provided that its boundary $\partial U$ is a smooth
one-codimensional submanifold of $M$.

\begin{cor}
Let $M$ be a separable Hilbert manifold. Then, for every two
disjoint closed subsets $C_1$, $C_2$ of $M$, there exists a
$C^\infty$ smooth function $\varphi:M\longrightarrow\mathbb{R}$
with no critical points, such that the level set
$N=\varphi^{-1}(0)$ is a $1$-codimensional $C^\infty$ smooth
submanifold of $M$ that separates $C_{1}$ and $C_{2}$, in the
following sense. Define $U_{1}=\{x\in M : \varphi(x)<0\}$ and
$U_{2}=\{x\in M : \varphi(x)>0\}$; then $U_1$ and $U_2$ are
disjoint $C^\infty$ smooth open sets of $M$ with common boundary
$\partial U_{1}=\partial U_{2}=N$, so that $C_{i}\subset U_{i}$
for $i=1, 2$.
\end{cor}
\begin{proof}
By Urysohn's lemma there exists a continuous function
$f:M\longrightarrow [0,1]$ so that $C_{1}\subset f^{-1}(0)$ and
$C_{2}\subset f^{-1}(1)$. Taking $\varepsilon=1/3$ and applying
Theorem \ref{main theorem for manifolds} we get a $C^\infty$
smooth function $\psi:M\longrightarrow\mathbb{R}$ which has no
critical points and is so that
    $$
    |f(x)-\psi(x)|\leq 1/3
    $$
for all $x\in M$; in particular
    $$
    C_{1}\subseteq f^{-1}(0)\subseteq \psi^{-1}(-\infty,1/2):=U_{1},
    $$
and
    $$
    C_{2}\subseteq f^{-1}(1)\subseteq \psi^{-1}(1/2, +\infty):=U_{2}.
    $$
The open sets $U_1$ and $U_2$ are smooth because their common
boundary $N=\psi^{-1}(1/2)$ is a smooth one-codimensional
submanifold of $M$ (thanks to the implicit function theorem and
the fact that $d\psi(x)\neq 0$ for all $x\in N$). In order to
obtain the result in the above form it is enough to set
$\varphi=\psi-1/2$.
\end{proof}

A trivial consequence of this result is that every closed subset
of $M$ can be uniformly approximated by smooth open subsets of
$M$. In fact,

\begin{cor}
Every closed subset of a separable Hilbert manifold $M$ can be
approximated by smooth open subsets of $M$, in the following
sense: for every closed set $C\subset M$ and every open set $W$
containing $C$ there is a $C^\infty$ smooth open set $U$ so that
$C\subset U\subseteq W$.
\end{cor}

Finally, the following result, which also implies the above
corollary, tells us that for every open set $U$ in $M$ there
always exists a function whose open support is $U$ and which does
not satisfy Rolle's theorem.

\begin{thm}\label{second main theorem}
For every open subset $U$ of a Hilbert manifold $M$ there is a
continuous function $f$ on $M$ whose support is the closure of
$U$, so that $f$ is $C^\infty$ smooth on $U$ and yet $f$ has no
critical point in $U$.
\end{thm}
\begin{proof}
For the same reasons as in the proof of Theorem \ref{main theorem
for manifolds} we may assume that $U$ is an open subset of the
Hilbert space $X=\ell_2$. Let $\varepsilon:X\longrightarrow
[0,+\infty)$ be the distance function to $X\setminus U$, that is,
    $$
    \varepsilon(x)=\textrm{dist}(x, X\setminus U)=\inf\{\|x-y\|: y\in X\setminus U\}.
    $$
The function $\varepsilon$ is continuous on $X$ and satisfies that
$\varepsilon(x)>0$ if and only if $x\in U$. According to Theorem
\ref{main theorem}, and setting $f(x)=2\varepsilon(x)$, there
exists a $C^\infty$ smooth function
$\psi:U\longrightarrow\mathbb{R}$ which has no critical points on
$U$, and such that $\varepsilon$-approximates $f$ on $U$, that is,
    $$
    |2\varepsilon(x)-\psi(x)|\leq\varepsilon(x)
    $$
for all $x\in U$. This inequality implies that
    $$
    \lim_{x\to z}\psi(x)=0
    $$
for every $z\in\partial U$. Therefore, if we set $\psi=0$ on
$X\setminus U$, the extended function
$\psi:X\longrightarrow[0,+\infty)$ is continuous on the whole of
$X$, is $C^\infty$ smooth on $U$ and has no critical points on
$U$. On the other hand, $\psi(x)\geq\varepsilon(x)>0$ for all
$x\in U$, hence the support of $\psi$ is $\overline{U}$.
\end{proof}

\section[Proofs]{Proof of the main result}

The main idea behind the proof of Theorem \ref{main theorem} is
as follows. First we use a {\em perturbed} smooth partition of unity to
approximate the given continuous mapping $f$. The summands of
this perturbed partition of unity are functions supported on
scalloped balls and carefully constructed in such a way that the
critical set $C_{\varphi}$ of the approximating sum $\varphi$
is locally compact.

Then we have to eliminate all the critical points without losing
much of the approximation. To this end we compose the
approximating mapping $\varphi$ with a deleting diffeomorphism
$h:X\longrightarrow X\setminus C_{\varphi}$ which extracts the
critical points $C_{\varphi}$ and is as close to the identity as
we want. The existence of such a diffeomorphism is guaranteed by
a quite elaborated result of West's \cite{West}. In this way we
obtain a smooth mapping $\psi$ which has no critical points, and
which happens to approximate the function $\varphi$ (which in
turn approximates the original $f$) because the perturbation
brought on $\varphi$ by the composition with $h$ is not very
important (recall that $h$ is arbitrarily closed to the identity).

The following proposition shows the existence of a function
$\varphi$ with the above properties. Recall that $C_{\varphi}$
stands for the set of critical points of $\varphi$.

\begin{prop}\label{existence of varphi}
Let $U$ be an open subset of the separable Hilbert space $X$. Let
$f:U\longrightarrow\mathbb{R}^{m}$ be a continuous mapping, and
$\varepsilon:U\longrightarrow (0,\infty)$ a continuous positive
function. Then there exist a $C^\infty$ smooth mapping
$\varphi:U\longrightarrow\mathbb{R}^{m}$ so that
\begin{itemize}
\item[{(a)}] $C_{\varphi}$ is locally compact and closed (relatively to U);

\item[{(b)}] $\|\varphi(x)-f(x)\|\leq \varepsilon(x)/2$ for all
$x\in U$.
\end{itemize}
\end{prop}

In fact, when the function $f$ takes values in the real line, we
can obtain a much stronger result which is interesting in itself
and might have some applications beyond the problem we are
dealing with, as it provides much more accurate information about
the structure and location of the critical points of the
approximation $\varphi$. The following theorem shows that any
continuous function can be uniformly approximated by $C^\infty$
smooth functions whose sets of critical points consist of
countable union of compact sets which are separated by pairwise
disjoint arbitrarily small open sets.
\begin{thm}\label{countable union of pairwise disjoint compact sets
of critical points} Let $U$ be an open subset of the separable
Hilbert space $X$. Let $f:U\longrightarrow\mathbb{R}$ be a
continuous function on $X$, and $\varepsilon:U\longrightarrow
(0,\infty)$ a continuous positive function. Then there exist a
$C^\infty$ smooth function $\varphi:U\longrightarrow\mathbb{R}$,
sequences $(K_{n})$ and $(U_{n})$ of compact sets and open sets
respectively, and a sequence $(B(y_{n}, r_{n}))$ of open balls
which are contained in $U$ and whose union covers $U$, such that:
\begin{itemize}
\item[{(a)}] $C_{\varphi}\subseteq\bigcup_{n=1}^{\infty} K_{n}$;
\item[{(b)}] $K_{n}\subset U_{n}\subseteq B(y_{n}, r_{n})$ for all
$n$, and $U_{n}\cap U_{m}=\emptyset$ whenever $n\neq m$;
\item[{(c)}] $|\varphi(x)-f(x)|\leq\varepsilon(x)$ for all $x$, and
$|\varphi(y)-f(x)|\leq\varepsilon(y_{n})$ for every $x,y\in
B(y_{n}, r_{n})$ and every $n$;
\item[{(d)}] for every $x\in U$ there exists an open neighborhood
$V_{x}$ of $x$ such that, either $V_{x}=U_{n}$ for a unique
$n=n_{x}$, or else $V_{x}\cap U_{n}=\emptyset$ for all $n$.
\end{itemize}
Moreover, for any given $r>0$, the radii of the balls can be
chosen so that $r_{n}\leq r$ for all $n$.
\end{thm}

Finally, the following restatement of a striking result of West's
\cite{West} ensures the existence of the diffeomorphism $h$. We
say that a mapping $g$ from a subset $A$ of $M$ is limited by an
open cover $G$ of $M$ if the collection $\{\{x,g(x)\}\, :\, x\in
A\}$ refines $G$.

\begin{thm}[West]\label{removing locally compact sets}
Let $C$ be a closed, locally compact subset of a Hilbert manifold
$M$, $U$ an open subset of $M$ with $C\subset U$, and $G$ an open
cover of $M$. Then there is a $C^\infty$ diffeomorphism $h$ of
$M$ onto $M\setminus C$ which is the identity outside $U$ and is
limited by $G$.
\end{thm}

Assume for a while that Proposition \ref{existence of varphi} is
already established, and let us see how we can deduce Theorem
\ref{main theorem}.

\medskip

\begin{center}
{\bf Proof of Theorem \ref{main theorem}}
\end{center}

For the given continuous mappings $f$ and $\varepsilon$, take a
mapping $\varphi$ with the properties of Proposition
\ref{existence of varphi}. Since $\varphi$ and $\varepsilon$ 
are continuous, for
every $z\in U$ there exists $\delta_{z}>0$ so that if $x,y\in
B(z,\delta_{z})$ then
$$\|\varphi(y)-\varphi(x)\|\leq\varepsilon(z)/4\leq\varepsilon(x)/2.$$

Let $G=\{B(x,\delta_{x})\, :\, x\in U\}$, $M=U$, and for the
critical set $C=C_{\varphi}$, use Theorem \ref{removing locally
compact sets} to find a $C^{\infty}$ diffeomorphism
$h:U\longrightarrow U\setminus C$ so that $h$ is limited by $G$.
Define $\psi=\varphi\circ h$.

Since $h$ is limited by $G$ we have that, for any given $x\in U$, there
exists $z\in U$ such that $x, h(x)\in B(z,\delta_{z})$, and therefore
$\|\varphi(h(x))-\varphi(x)\|\leq\varepsilon(z)/4,$ that is, we
have that
      $$
      \|\psi(x)-\varphi(x)\|\leq\varepsilon(z)/4\leq\varepsilon(x)/2.
      $$
Hence, by combining this inequality with (b) of
Proposition \ref{existence of varphi}, we obtain that $$
\|\psi(x)-f(x)\|\leq\varepsilon(x) \eqno(1) $$ for all $x\in U$.

Let us see that $\psi$ does not have any
critical point.
The derivative of $\psi$ is given by
    $$
    d\psi(x)=d\varphi(h(x))\circ d h(x). \eqno(2)
    $$
Since $h(x)\notin C=C_{\varphi}$, we have that the linear map
$d\varphi(h(x))$ is surjective. On the other hand $d h(x)$ is a
linear isomorphism (because $h$ is a diffeomorphism). Then it is
clear that the composition $d\psi(x)=d\varphi(h(x))\circ d h(x)$
is a linear surjection from $X$ onto $\mathbb{R}^{m}$, for every
$x\in U$. \hspace{9cm}

\begin{rem}
{\em In the case when $f:U\longrightarrow\mathbb{R}$ we do not
need to use the full power of West's result. Thanks to the more
accurate statement provided by Theorem \ref{countable union of
pairwise disjoint compact sets of critical points} we can instead
use a much more elementary result that tells us that for every
compact subset $K$ and every open subset $U$ of $X$ with $K\subset
U$, there exists a $C^\infty$ diffeomorphism $h:X\longrightarrow
X\setminus K$ such that $h$ restricts to the identity outside $U$.
In our case, to eliminate the critical points of the
approximating function $\varphi$ of Theorem \ref{countable union
of pairwise disjoint compact sets of critical points}, we may
compose $\varphi$ with a sequence of deleting diffeomorphisms
$h_{n}:X\longrightarrow X\setminus K_{n}$ which extract each of
the compact sets of critical points $K_n$ and restrict to the
identity outside each of the open sets $U_{n}$. The infinite
composition of deleting diffeomorphisms with our function,
$\psi=\varphi\circ \bigcirc_{n=1}^{\infty} h_{n}$, is locally
finite, in the sense that only a finite number (in fact at most
one) of the diffeomorphisms are acting on some neighborhood of
each point, while all the rest restrict to the identity on that
neighborhood. As in the proof above, it follows that $\psi$ has no
critical points (we can use exactly the same argument locally), 
and still approximates $f$ (recall that each
$h_n$ restricts to the identity outside the set $U_n$, on which
$\varphi$ has a very small oscillation, and the $U_n$ are pairwise
disjoint).}
\end{rem}

\medskip

\begin{center}
{\bf Proof of Proposition \ref{existence of varphi}}
\end{center}

We will assume that $U=X$, since the proof is completely analogous
in the case of a general open set. One only has to take some (easy
but rather rambling) technical precautions in order to make sure
that the different balls considered in the argument are in $U$.

In order to avoid bearing an unnecessary burden of notation, we
will make the proof of this proposition for the case of a constant
$\varepsilon>0$. Later on we will briefly explain what additional
technical precautions must be taken in order to deduce the
general form of this result (see Remark \ref{remark for e(x)}
below).

Let $B(x, r)$ and $\overline{B}(x,r)$ stand for the open ball and
closed ball, respectively, of center $x$ and radius $r$, with
respect to the usual hilbertian norm $\|\cdot\|$ of $X$.

\medskip

\noindent {\bf Case I.} We will first consider the case of a real
valued function $f:U\longrightarrow\mathbb{R}$. Fix
$\varepsilon>0$. By continuity, for every $x\in X$ there exists
$\delta_{x}>0$ so that $|f(y)-f(x)|\leq \varepsilon/8$ whenever
$y\in B(x, 2\delta_{x})$. Since $X=\bigcup_{x\in X}
B(x,\delta_{x}/2)$ is separable, there exists a countable
subcovering,
    $$
    X=\bigcup_{n=1}^{\infty}B(x_{n}, r_{n}/2),
    $$
where $r_{n}=\delta_{x_{n}}$, for some sequence of centers
$(x_{n})$. By induction (and using the fact that every
finite-dimensional subspace of $X$ has empty interior in $X$), we
can choose a sequence of {\em linearly independent} vectors
$(y_{n})$, with $y_{n}\in B(x_{n}, r_{n}/2)$, so that
    $$
    X=\bigcup_{n=1}^{\infty}B(y_{n}, r_{n}). \eqno(3)
    $$
Moreover, we have that
    $$
    |f(y)-f(y_{n})|\leq\varepsilon/4\,
    \text{ whenever $\|y-y_{n}\|\leq r_{n}$.} \eqno (4)
    $$
Now we define the scalloped balls $B_{n}$ that are the basis for
our perturbed partition of unity: set $B_{1}=B(y_{1},r_{1})$, and
for $n\geq 2$ define
    $$
    B_{n}=B(y_{n}, r_{n})\setminus \Bigl(\bigcup_{j=1}^{n-1}
    \overline{B}(y_{j}, \lambda_{n} r_{j})\Bigr);
    $$
where $1/2<\lambda_{2}<\lambda_{3}< ...
<\lambda_{n}<\lambda_{n+1}< ... <1$, with $\lim_{n\to\infty}\lambda_{n}=1$.

Taking into account that $\lim_{n\to\infty}\lambda_{n}=1$, it is
easily checked that the $B_{n}$ form a locally finite open
covering of $X$, with the nice property that
$$|f(y)-f(y_{n})|\leq\varepsilon/4 \, \text{ whenever } y\in
B_{n}.$$

Next, pick a $C^{\infty}$ smooth function
$g_{1}:\mathbb{R}\longrightarrow [0,1]$ so that:
\begin{itemize}
\item[{(i)}] $g_{1}(t)=1$ for $t\leq 0$,
\item[{(ii)}] $g_{1}(t)=0$ for $t\geq {r_{1}}^{2}$,
\item[{(iii)}] $g_{1}'(t)<0$ if $0<t< {r_{1}}^{2}$;
\end{itemize}
and define then $\varphi_{1}:X\longrightarrow\mathbb{R}$ by
    $$
    \varphi_{1}(x)=g_{1}(\|x-y_{1}\|^{2})
    $$
for all $x\in X$. Note that $\varphi_{1}$ is a $C^{\infty}$ smooth
function whose open support is $B_{1}$, and $B_{1}\cap
C_{\varphi_{1}}=\{y_{1}\}$, that is, $y_{1}$ is the only critical
point of $\varphi_{1}$ that lies inside $B_{1}$.

Now, for $n\geq 2$, pick $C^\infty$ smooth functions
$\theta_{(n,j)}:\mathbb{R}\longrightarrow [0,1]$, $j=1, ..., n$,
with the following properties. For $j=1, ..., n-1$,
$\theta_{(n,j)}$ satisfies that
\begin{itemize}
\item[{(i)}] $\theta_{(n,j)}(t)=0$ for $t\leq(\lambda_{n}r_{j})^{2}$,
\item[{(ii)}] $\theta_{(n,j)}(t)=1$ for $t\geq {r_{j}}^{2}$,
\item[{(iii)}] $\theta_{(n,j)}'(t)>0$ if $(\lambda_{n}r_{j})^{2}<t<{r_{j}}^{2}$;
\end{itemize}
while for $j=n$ the function $\theta_{(n,n)}$ is such that
\begin{itemize}
\item[{(i)}] $\theta_{(n,n)}(t)=1$ for $t\leq 0$,
\item[{(ii)}] $\theta_{(n,n)}(t)=0$ for $t\geq {r_{n}}^{2}$,
\item[{(iii)}] $\theta_{(n,n)}'(t)<0$ if $0<t< {r_{n}}^{2}$.
\end{itemize}
Then define the function $g_{n}:\mathbb{R}^{n}\longrightarrow
[0,1]$ as
    $$
    g_{n}(t_{1}, ..., t_{n})=\prod_{i=1}^{n}\theta_{(n,i)}(t_{i})
    $$
for all $t=(t_{1}, ..., t_{n})\in\mathbb{R}^{n}$. This function
is clearly $C^\infty$ smooth on $\mathbb{R}^{n}$ and satisfies
the following properties:
\begin{itemize}
\item[{(i)}] $g_{n}(t_{1},...,t_{n})>0$ if and only if
$t_{j}>(\lambda_{n}r_{j})^{2}$ for all $j=1, ..., n-1$, and
$t_{n}< {r_{n}}^{2}$; and $g_{n}$ vanishes elsewhere;
\item[{(ii)}] $g_{n}(t_{1},...,t_{n})=\theta_{(n,n)}(t_{n})$ whenever $t_{j}\geq
{r_{j}}^{2}$ for all $j=1, ..., n-1$;
\item[{(iii)}] $\nabla g_{n}(t_{1},...,t_{n})\neq 0$ provided
$(\lambda_{n}r_{j})^{2}<t_{j}$ for all $j=1, ..., n-1$, and
$0<t_{n}<{r_{n}}^{2}$.
\end{itemize}
Moreover, under the same conditions as in (iii) just above we
have that
    $$
    \frac{\partial g_{n}}{\partial t_{n}}(t_{1},...,t_{n})=
    \frac{\partial \theta_{(n,n)}}{\partial t_{n}}(t_{n})
    \prod_{i=1}^{n-1}\theta_{(n,i)}(t_{i})<0, \eqno(5)
    $$
since no function in this product vanishes on the specified set,
while for $j<n$, according to the corresponding properties of the
functions $\theta_{(n,j)}$ we have that
    $$
    \frac{\partial g_{n}}{\partial t_{j}}(t_{1},...,t_{n})=
    \frac{\partial \theta_{(n,j)}}{\partial t_{j}}(t_{j})
    \prod_{i=1, i\neq j}^{n}\theta_{(n,i)}(t_{i})>0. \eqno(6)
    $$
If we are not in the conditions of (iii) then the corresponding
inequalities do still hold but are not strict.

Let us now define $\varphi_{n}:X\longrightarrow [0,1]$ by
    $$
    \varphi_{n}(x)=g_{n}(\|x-y_{1}\|^{2}, ... , \|x-y_{n}\|^{2}).
    $$
It is clear that $\varphi_{n}$ is a $C^{\infty}$ smooth function
whose open support is precisely the scalloped ball $B_{n}$.

As above, let us denote by $C_{\varphi_{n}}$ the critical set of
$\varphi_{n}$, that is,
$$
C_{\varphi_{n}}=\{x\in X :
\varphi_{n}'(x)=0\}.
$$
Since our norm $\|\cdot\|$ is hilbertian we
have that, if $x\in C_{\varphi_{n}}\cap B_{n}$, then $x$ belongs
to the affine span of $y_{1}, ..., y_{n}$. Indeed, if $x\in
B_{n}$,
    $$
    \varphi_{n}'(x)=\sum_{j=1}^{n}\frac{\partial g_{n}}{\partial t_{j}}
    (\|x-y_{1}\|^{2}, ... , \|x-y_{n}\|^{2})\,2(x-y_{j})=0, \eqno(7)
    $$
which (taking into account (5) and the fact that the $y_{j}$ are
all linearly independent) means that $x$ is in the affine span of
$y_{1}, ..., y_{n}$. Here, as is usual, we identify the Hilbert
space $X$ with its dual $X^*$, and we make use of the fact that
the derivative of the function $x\mapsto \|x\|^{2}$ is the
mapping $x\mapsto 2x$.

Similarly, by using (5) it can be shown that $x\in
C_{\varphi_{1}+\cdots+\varphi_{m}}\cap (B_{1}\cup ...\cup B_{m})$
implies that $x$ belongs to the affine span of $y_{1}, ...,
y_{m}$.

\medskip

In order that our approximating function has a small critical set
we cannot use the standard approximation provided by the partition
of unity associated with the functions
$(\varphi_{j})_{i\in\mathbb{N}}$, namely
    $$
    x\mapsto \frac{\sum_{n=1}^{\infty}\alpha_{n}\varphi_{n}(x)}
    {\sum_{n=1}^{\infty}\varphi_{n}(x)},
    $$
where $\alpha_{n}=f(y_{n})$. Indeed, such a function would have a
huge set of critical points since it would be constant (equal to
$\alpha_{n}$) on a lot of large places (at least on each $B_{n}$
minus the union of the rest of the $B_j$). Instead, we will
modify this standard approximation by letting the $\alpha_{n}$ be
functions (and not mere numbers) of very small oscillation and
with only one critical point (namely $y_{n}$). So, for every
$n\in\mathbb{N}$ let us pick a $C^{\infty}$ smooth real function
$a_{n}:[0, +\infty)\longrightarrow\mathbb{R}$ with the following
properties:
\begin{itemize}
\item[{(i)}] $a_{n}(0)=f(y_{n})$;
\item[{(ii)}] $a_{n}'(t)<0$ whenever $t>0$;
\item[{(iii)}] $|a_{n}(t)-a_{n}(0)|\leq\varepsilon/4$ for all $t\geq 0$;
\end{itemize}
and define $\alpha_{n}:X\longrightarrow\mathbb{R}$ by
    $$
    \alpha_{n}(x)=a_{n}(\|x-y_{n}\|^{2})
    $$
for every $x\in X$. It is clear that $\alpha_{n}$ is a $C^\infty$
smooth function on $X$ whose only critical point is $y_{n}$.
Besides, $$ |\alpha_{n}(x)-f(y_{n})|\leq\varepsilon/4 \text{ for
all $x\in X$}. $$

Now we can define our approximating function
$\varphi:X\longrightarrow\mathbb{R}$ by
    $$
    \varphi(x)=\frac{\sum_{n=1}^{\infty}\alpha_{n}(x)\varphi_{n}(x)}
    {\sum_{n=1}^{\infty}\varphi_{n}(x)}
    $$
for every $x\in X$. Since the sums are locally finite, it is clear
that $\varphi$ is a well-defined $C^\infty$ smooth function.

\begin{fact}\label{good approximation}
The function $\varphi$ approximates $f$ nicely. Namely, we have
that
\begin{itemize}
\item[{(i)}] $|\varphi(x)-f(x)|\leq\varepsilon/2$ for all $x\in X$, and
\item[{(ii)}] $|\varphi(y)-f(x)|\leq \varepsilon$ for all $x, y\in
B(y_{n}, r_{n})$ and each $n\in\mathbb{N}$.
\end{itemize}
\end{fact}
\begin{proof}
For every $n$ we have that
$|\alpha_{n}(x)-f(y_{n})|\leq\varepsilon/4$ for all $x\in X$. On
the other hand, by $(4)$ above we know that $|f(x)-f(y_{n})|\leq
\varepsilon/4$ whenever $x\in B(y_{n}, r_{n})$. Then, by the
triangle inequality, it follows that
    $$
    |\alpha_{n}(x)-f(x)|\leq\varepsilon/2 \eqno(8)
    $$
whenever $x\in B(y_{n}, r_{n})$. In the same way we deduce that
    $$
    |\alpha_{m}(x)-f(y_{n})|\leq\varepsilon/2 \eqno(9)
    $$
whenever $x\in B(y_{n}, r_{n})\cap B(y_{m}, r_{m})$. Since
$\varphi_{m}(y)=0$ when $y\notin B(y_{m}, r_{m})$, from $(8)$ we
get that
    $$
    |\varphi(x)-f(x)|=
    \bigg|\frac{\sum_{m=1}^{\infty}(\alpha_{m}(x)-f(x))\varphi_{m}(x)}
    {\sum_{m=1}^{\infty}\varphi_{m}(x)}\bigg|\leq
    \frac{\sum_{m=1}^{\infty}\frac{\varepsilon}{2} \varphi_{m}(x)}
    {\sum_{m=1}^{\infty}\varphi_{m}(x)}=\varepsilon/2
    $$
for all $x\in X$, which shows $(i)$. Similarly, we deduce from
$(9)$ that
    $$
    |\varphi(y)-f(y_{n})|=
    \bigg|\frac{\sum_{m=1}^{\infty}(\alpha_{m}(y)-f(y_{n}))\varphi_{m}(y)}
    {\sum_{m=1}^{\infty}\varphi_{m}(y)}\bigg|\leq
    \frac{\sum_{m=1}^{\infty}\frac{\varepsilon}{2} \varphi_{m}(y)}
    {\sum_{m=1}^{\infty}\varphi_{m}(y)}=\varepsilon/2
    $$
for every $y\in B(y_{n}, r_{n})$, which, combined with $(4)$
above, yields that
    $$
    |\varphi(y)-f(x)|\leq \varepsilon/2+\varepsilon/4,
    $$
for every $x, y\in B(y_{n}, r_{n})$, so $(ii)$ is satisfied as
well.
\end{proof}

\medskip

Now let us have a look at the derivative of $\varphi$. To this
end let us introduce the auxiliary functions $f_{n}$ defined by
    $$
    f_{n}(x)=\frac{\sum_{k=1}^{n}\alpha_{k}(x)\varphi_{k}(x)}
    {\sum_{k=1}^{n}\varphi_{k}(x)}, \, \textrm{ for all }\,
    x\in \bigcup_{i=1}^{n}B_{i}.
    $$
Notice that $\varphi$ can be expressed as
    $$
    \varphi(x)=\lim_{n\to\infty}f_{n}(x),
    $$
that the domains of the $f_{n}$ form an increasing tower of open
sets whose union is $X$, and that each $f_{n}$ restricts to
$f_{n-1}$ on $\bigcup_{i=1}^{n-1}B_{i}\setminus B_{n}$. Moreover,
we have the following.
\begin{fact}\label{varphi is locally f_n}
For each $x\in X$ there is an open neighborhood $V_{x}$ of $x$ and
some $n_{x}\in\mathbb{N}$ so that $\varphi(y)=f_{n}(y)$ for all
$y\in V_{x}$ and all $n\geq n_{x}$.
\end{fact}
\begin{proof}
Indeed, we have that, for every $n\in\mathbb{N}$,
    $$
    \varphi(y)=f_{k}(y) \hspace{0.3cm} \textrm{for all} \hspace{0.3cm}
    y\in V_{n}:=\big(\bigcup_{j=1}^{n}B_{j}\big)\setminus
    \big(\bigcup_{i=n+1}^{\infty}\overline{B}_{i}\big),
    \hspace{0.3cm} \textrm{and for all} \hspace{0.3cm} k\geq n.
    $$
The $V_{n}$ are open, $V_{n}\subseteq V_{n+1}$, and
$\bigcup_{i=1}^{\infty}V_{i}=X$, because the covering of $X$
formed by the $B_j$ is locally finite.
\end{proof}

Hence, by looking at the derivatives of the functions $f_{n}$ we
will get enough information about the derivative of $\varphi$.

If $x\in\bigcup_{j=1}^{n}B_{j}$ then the expression for the
derivative of $f_{n}$ is given by
    $$
    f_{n}'(x)=\frac{\sum_{j=1}^{n}[\alpha_{j}'(x)\varphi_{j}(x)+
    \alpha_{j}(x)\varphi_{j}'(x)]\sum_{i=1}^{n}\varphi_{i}(x)-
    \sum_{j=1}^{n}\varphi_{j}'(x)\sum_{i=1}^{n}\alpha_{i}(x)\varphi_{i}(x)}
    {(\sum_{j=1}^{n}\varphi_{j}(x))^{2}}.
    $$
Therefore, for $x\in\bigcup_{j=1}^{n}B_{j}$ we have that
$f_{n}'(x)=0$ if and only if
    $$
    \sum_{j=1}^{n}\sum_{i=1}^{n}\varphi_{i}(x)
    \Bigl[\alpha_{j}'(x)\varphi_{j}(x)+\bigl(\alpha_{j}(x)-\alpha_{i}(x)\bigr)\varphi_{j}'(x)\Bigr]
    =0. \eqno(10)
    $$
By inserting the expressions for the derivatives of $\varphi_{j}$
and $\alpha_{j}$ in equation $(10)$, we can express the condition
$f_{n}'(x)=0$ as a nontrivial linear dependence link on the
vectors $(x-y_{j})$, which yields that $x$ is in the affine span
of the points $y_{1}, ..., y_{n}$.
\begin{notation}
{\em In the sequel $\mathcal{A}[z_{1},...,z_{k}]$ stands for the
affine subspace spanned by a finite sequence of points $z_{1},
..., z_{k}\in X$.}
\end{notation}
\begin{fact}\label{the critical set is locally finite
dimensional} If $x\in C_{f_{n}}\cap B_{n}$, then $x\in
\mathcal{A}_{n}:=\mathcal{A}[y_{1}, ..., y_{n}]$. Moreover, for
each $n\in\mathbb{N}$ and for every finite sequence of positive
integers $k_{1}<k_{2}<...<k_{m}<n$ we have that
    $$
    C_{f_{n}}\cap\big(B_{n}\setminus\bigcup_{j=1}^{m}B_{k_{j}}\big)
    \subseteq\mathcal{A}\big[\{y_{1}, ..., y_{n}\}
    \setminus\{y_{k_{1}}, ..., y_{k_{m}}\}\big].
    $$
\end{fact}
\begin{proof}
As above, in all the subsequent calculations, we will identify the
Hilbert space $X$ with its dual $X^*$, and the derivative of
$\|\cdot\|^{2}$ with the mapping $x\mapsto 2x$. To save notation,
let us simply write
    $$
    \frac{\partial g_{n}}{\partial
    t_{j}}(\|x-y_{1}\|^{2},...,\|x-y_{n}\|^{2})=\mu_{(n,j)},
    $$
and
    $$
    a_{j}'(\|x-y_{j}\|^{2})=\eta_{j}.
    $$

Notice that, according to $(5)$ and $(6)$ above, $\mu_{(n,j)}\geq
0$ for $j=1, ..., n-1$, while $\mu_{(n,n)}\leq 0$; and
$\mu_{(n,n)}\neq 0$ provided $x\in B_{n}$ and $x\neq y_{n}$; on
the other hand it is clear that $\eta_{j}<0$ for all $j$ unless
$x=y_{j}$ (in which case $\eta_{j}=0$).

Assuming $x\in C_{f_{n}}\cap B_{n}$, and taking into account the
expression $(10)$ for $\varphi_{j}'(x)$ and the fact that
$\alpha_{j}'(x)=2\eta_{j}(x-y_{j})$, we can write condition $(10)$
above in the form
    $$
    2\sum_{j=1}^{n}\sum_{i=1}^{n}\varphi_{i}(x)
    \Bigl[\eta_{j}\varphi_{j}(x)\,(x-y_{j})+\bigl(\alpha_{j}(x)-\alpha_{i}(x)\bigr)
    \sum_{\ell=1}^{j}\mu_{(j,\ell)}\,(x-y_{\ell})\Bigr]
    =0,
    $$
which in turn is equivalent (taking the common factors of each
$(x-y_{j})$ together) to the following one
    $$
    \sum_{j=1}^{n}\Biggl[
    \eta_{j}\varphi_{j}(x)\sum_{i=1}^{n}\varphi_{i}(x) +
    \sum_{k=j}^{n}\Bigl(\sum_{i=1}^{n}\bigl(\alpha_{k}(x)-
    \alpha_{i}(x)\bigr)\varphi_{i}(x)\Bigr)
    \mu_{(k,j)}\Biggr]\,(x-y_{j})=0. \eqno(11)
    $$
Now notice that, if we can prove that at least one of the
expressions multiplying the $(x-y_{j})$ does not vanish then we
are done; indeed, we will have that the vectors $x-y_{1}$, ...,
$x-y_{n}$ are linearly dependent, which means that $x$ belongs to
the affine span of the points $y_{1}, ..., y_{n}$.

So let us check that not all of those expressions in $(11)$
vanish. In fact we are going to see that at least one of the
terms is strictly negative. We can obviously assume that $x$ is
not any of the points $y_{1}, ..., y_{n}$ (which are already in
$\mathcal{A}_{n}$). In this case we have that $\mu_{(n,n)}<0$ and
$\eta_{j}<0$ for all $j=1, ..., n$. For simplicity, we will only
make the argument in the case $n=3$; giving a proof in a more
general case would be as little instructive as tedious to read.

Let us first assume that $\varphi_{j}(x)\neq 0$ for $j=1,2,3$. We
begin by looking at the term that multiplies $(x-y_{3})$ in
$(11)$, that is
    $$
    \beta_{3}:=\eta_{3}\varphi_{3}(x)\sum_{i=1}^{3}\varphi_{i}(x) +
    \sum_{i=1}^{3}\bigl(\alpha_{3}(x)-
    \alpha_{i}(x)\bigr)\varphi_{i}(x)
    \mu_{(3,3)}.
    $$
If $\sum_{i=1}^{3}\bigl(\alpha_{3}(x)-\alpha_{i}(x)\bigl)
\varphi_{i}(x)\geq 0$ we are done, since in this case we easily
see that $\beta_{3}<0$ (remember that $\mu_{(3,3)}\leq 0$,
$\eta_{3}<0$, and $\varphi_{3}(x)>0$). Otherwise we have that $$
\sum_{i=1}^{3}\Bigl(\alpha_{3}(x)-\alpha_{i}(x)\Bigr)
\varphi_{i}(x)<0, $$ and then we look at the term $\beta_{2}$
multiplying $(x-y_{2})$ in $(11)$, namely,
    $$
    \beta_{2}:=\eta_{2}\varphi_{2}(x)\sum_{i=1}^{3}\varphi_{i}(x) +
    \sum_{k=2}^{3}\Biggl(\sum_{i=1}^{3}\Bigl(\alpha_{k}(x)-
    \alpha_{i}(x)\Bigr)\varphi_{i}(x)\Biggr)
    \mu_{(k,2)}.
    $$
Now, since $\mu_{(3,2)}\geq 0$, we have
$\sum_{i=1}^{3}\bigl(\alpha_{3}(x)-\alpha_{i}(x)\bigr)
\varphi_{i}(x)\mu_{(3,2)}\leq 0$, and on the other hand
$\eta_{2}\varphi_{2}(x)\sum_{i=1}^{3}\varphi_{i}(x)<0$ so that,
if $\sum_{i=1}^{3}\bigl(\alpha_{2}(x)-\alpha_{i}(x)\bigr)
\varphi_{i}(x)$ happens to be nonnegative, then we also have
$\sum_{i=1}^{3}\bigl(\alpha_{2}(x)-\alpha_{i}(x)\bigr)
\varphi_{i}(x)\mu_{(2,2)}\leq 0$, and then we are done since
$\beta_{2}$, being a sum of negative terms (one of them strictly
negative) must be negative as well. Otherwise, $$
\sum_{i=1}^{3}\Bigl(\alpha_{2}(x)-\alpha_{i}(x)\Bigr)
\varphi_{i}(x) $$ is negative, and then we finally pass to the
term $\beta_{1}$ multiplying $(x-y_{1})$ in $(11)$, that is,
    $$
    \beta_{1}:=\eta_{1}\varphi_{1}(x)\sum_{i=1}^{3}\varphi_{i}(x) +
    \sum_{k=1}^{3}\Biggl(\sum_{i=1}^{3}\Bigl(\alpha_{k}(x)-
    \alpha_{i}(x)\Bigr)\varphi_{i}(x)\Biggr)
    \mu_{(k,1)}.
    $$
Here, by the assumptions we have made so far and taking into
account the signs of $\mu_{(k,j)}$ and $\eta_{j}$, we see that
$\sum_{i=1}^{3}\bigl(\alpha_{k}(x)- \alpha_{i}(x)\bigr)\varphi_{i}(x)
\mu_{(k,1)}\leq 0$ for $k=2, 3$. Having arrived at this point, it
is sure that
$\sum_{i=1}^{3}\bigl(\alpha_{1}(x)-\alpha_{i}(x)\bigr)\varphi_{i}(x)$ must
be nonnegative (otherwise the numbers
$\sum_{i=1}^{3}\bigl(\alpha_{k}(x)-\alpha_{i}(x)\bigr)\varphi_{i}(x)$
should be strictly negative for all $k=1,2,3$, which is
impossible if one takes $\alpha_{k}(x)$ to be the maximum of the
$\alpha_{i}(x)$), and now we can deduce as before that
$\beta_{1}<0$.

Finally let us consider the case when some of the $\varphi_{i}(x)$
vanish, for $i=1,2$ (remember that $\varphi_{3}(x)\neq 0$ since
$x\in B_{3}$, the open support of $\varphi_{3}$). From the
definitions of $\mu_{(k,j)}$, $g_{n}$ and $\varphi_{n}$, it is
clear that $\mu_{(k,j)}=0$ whenever $\varphi_{j}(x)=0$ or
$\varphi_{k}(x)=0$, and bearing this fact in mind we can simplify
equality $(11)$ to a great extent by dropping all the terms that
now vanish.

If $\varphi_{1}(x)=\varphi_{2}(x)=0$ then $(11)$ reads
    $$
    \varphi_{3}(x)^{2}\eta_{3}\,(x-y_{3})=0,
    $$
which cannot happen since we assumed $x\neq y_{j}$ (this means
that the only critical point that $f_{n}$ can have in
$B_{3}\setminus(B_{1}\cup B_{2})$ is $y_{3}$).

If $\varphi_{1}(x)=0$ and $\varphi_{2}(x)\neq 0$ then the term
$\beta_{1}$ accompanying $(x-y_{1})$ in $(11)$ vanishes, and hence
$(11)$ is reduced to
    $$
    \sum_{j=2}^{3}\Biggl[
    \eta_{j}\varphi_{j}(x)\sum_{i=2}^{3}\varphi_{i}(x) +
    \sum_{k=j}^{3}\Bigl(\sum_{i=2}^{3}(\alpha_{k}(x)-
    \alpha_{i}(x))\varphi_{i}(x)\Bigr)
    \mu_{(k,j)}\Biggr]\,(x-y_{j})=0.
    $$
Since at least one of the numbers $\sum_{i=2}^{3}(\alpha_{k}(x)-
\alpha_{i}(x))\varphi_{i}(x)$, $k=2,3$, is nonnegative, the same
reasoning as in the first case allows us to conclude that either
$\beta_{3}$ or $\beta_{2}$ is strictly negative. Finally, in the
case $\varphi_{1}(x)\neq 0$ and $\varphi_{2}(x)=0$, it is
$\beta_{2}$ that vanishes, and $(11)$ reads $\beta_{1}\,(x-y_{1})+
\beta_{3}\,(x-y_{3})=0$, where
    $$
    \beta_{3}=\eta_{3}\varphi_{3}(x)\sum_{i=1, i\neq 2}^{3}\varphi_{i}(x) +
    \sum_{i=1, i\neq 2}^{3}\Bigl(\alpha_{3}(x)-
    \alpha_{i}(x)\Bigr)\varphi_{i}(x)
    \mu_{(3,3)},
    $$
and $$
    \beta_{1}=\eta_{1}\varphi_{1}(x)\sum_{i=1, i\neq 2}^{3}\varphi_{i}(x) +
    \sum_{k=1, i\neq 2}^{3}\sum_{i=1, i\neq 2}^{3}\Bigl(\alpha_{k}(x)-
    \alpha_{i}(x)\Bigr)\varphi_{i}(x)
    \mu_{(k,1)}.
    $$
Again, at least one of the numbers $\sum_{i=1, i\neq
2}^{3}\bigl(\alpha_{k}(x)- \alpha_{i}(x)\bigr)\varphi_{i}(x)$, $k=1,3$, is
nonnegative, and the same argument as above applies.

Finally, bearing in mind the definition of the functions
$\varphi_{j}$, whose open support are the $B_{j}$, it is clear
that the above discussion shows, in fact, the following
inclusions:
\begin{itemize}
\item[{}] $C_{f_{3}}\cap B_{3}\subseteq\mathcal{A}[y_{1}, y_{2},
y_{3}]$;
\item[{}] $C_{f_{3}}\cap(B_{3}\setminus B_{1})\subseteq\mathcal{A}[y_{2},
y_{3}]$, and $C_{f_{3}}\cap(B_{3}\setminus
B_{2})\subseteq\mathcal{A}[y_{1}, y_{3}]$ ;
\item[{}] $C_{f_{3}}\cap(B_{3}\setminus (B_{1}\cup
B_{2}))\subseteq \mathcal{A}[y_{3}]$.
\end{itemize}
An analogous argument in the case $n\geq 4$ proves the second part
of the statement of Fact~\ref{the critical set is locally finite
dimensional}.
\end{proof}
\begin{rem}\label{the derivative is a linear combination of x-y_n}
{\em Notice that the above proof shows that the derivative
$df_{n}(x)$ of the function $f_{n}$ at a point $x$ can be
expressed as a nontrivial linear combination of the linear
functionals $(x-y_{k})\in\ell_{2}^{*}=\ell_{2}$, $k=1,...,n$. That
is, for every $x\in\bigcup_{i=1}^{n}B_{i}$ there are numbers
$\beta_{1}(x), ..., \beta_{n}(x)$ such that at least one of them
does not vanish, and
    $$
    df_{n}(x)=\frac{1}{(\sum_{j=1}^{n}\varphi_{j}(x))^{2}}
    \sum_{k=1}^{n}\beta_{k}(x) (x-y_{k}).
    $$
This will turn out to be a crucial observation when dealing with
the case $m\geq 1$.}
\end{rem}

\medskip

Since $\varphi$ has a continuous derivative, it is obvious that
its critical set $C_{\varphi}$ is closed in $U$. According to Fact
\ref{varphi is locally f_n}, $\varphi$ locally coincides with one
of the $f_{n}$. From Fact~\ref{the critical set is locally finite
dimensional} it follows that the set of critical points
$C_{f_{n}}$ of each function $f_n$ is contained in a finite
dimensional affine subspace of $X$. Therefore it is clear that
the set $C_{\varphi}$ of critical points of $\varphi$ is locally
contained in finite dimensional subspaces, that is, for each $x\in
C_{\varphi}$ there is an open bounded neighborhood $V_{x}$ of $x$
so that $C_{\varphi}\cap\overline{V_{x}}$ is contained in a
finite dimensional subspace $F_{x}$ of $X$ and hence is compact
(as is closed and bounded as well). This means that $C_{\varphi}$
is locally compact, and concludes the proof of Proposition
\ref{existence of varphi} in the case $m=1$.

\medskip

\noindent {\bf Case II.} Let us now deal with the case when
$f:X\longrightarrow\mathbb{R}^{m}$ with $m\geq 2$. We denote
$f=(f^{1}, ..., f^{m})$, where $f^{1}, ..., f^{m}$ are the
coordinate functions of $f$. In this case we have to construct
$C^\infty$ smooth functions $\varphi^{1}, ..., \varphi^{m}$ so
that each $\varphi^{j}$ uniformly approximates $f^j$ and the set
of points $x\in X$ at which the derivatives $d\varphi^{1}(x)$,
..., $d\varphi^{m}(x)$ are linearly dependent is locally compact.
If we succeed in doing so then it is clear that the function
$\varphi=(\varphi^{1}, ...,
\varphi^{m}):X\longrightarrow\mathbb{R}^{m}$ will approximate $f$
and its set $C_{\varphi}$ of critical points will be closed and
locally compact.

Let us define $\varepsilon_{j}=\varepsilon/\sqrt{4m}$, $j=1, ...,
m$. As each of the functions $f^{j}$, with $j=1, ..., m$, is
continuous, for every $x\in X$ there exists $\delta^{j}_{x}>0$ so
that $$|f^{j}(y)-f^{j}(x)|\leq \varepsilon_{j}/8 \hspace{0.3cm}
\text{ for all } \hspace{0.3cm} y\in B(x, 2\delta_{x}).$$ Since
$X=\bigcup_{x\in X} B(x,\delta^{j}_{x}/2)$ is separable, we may
take a countable subcovering,
    $$
    X=\bigcup_{n=1}^{\infty}B(x^{j}_{n}, r^{j}_{n}/2),
    $$
where $r^{j}_{n}=\delta^{j}_{x_{n}}$, for each $j=1, ..., m$.

Now, we can slightly perturb the centers $x^{j}_{n}$ of the balls
so that the union of all the $m$ sequences of centers forms a set
of linearly independent vectors. Indeed, bearing in mind that the
complement of every finite dimensional subspace of $X$ is dense
in the infinite dimensional space $X$, we may inductively choose
(taking $m$ points $y^{1}_{k}, ..., y^{m}_{k}$ at each $k$-th step
of the induction process) sequences of points
$(y^{j}_{n})_{n=1}^{\infty}$, $j=1, ..., m$, with $y^{j}_{n}\in
B(x^{j}_{n}, r^{j}_{n}/2)$, so that:
\begin{enumerate}
\item[{(i)}] $\{y^{j}_{n} \, : \, n\in\mathbb{N}, \, j=1, ..., m\}$ is a
set of linearly independent vectors;
\item[{(ii)}]  $X=\bigcup_{n=1}^{\infty}B(y^{j}_{n}, r^{j}_{n})$ for every $j=1, ...,
m$; and
\item[{(iii)}]  $|f^{j}(y)-f^{j}(y^{j}_{n})|\leq\varepsilon_{j}/4$ whenever
$\|y-y^{j}_{n}\|\leq r^{j}_{n}$.
\end{enumerate}

Next, for each collection of balls $\{B(y^{j}_{n},
r^{j}_{n})\}_{n\in\mathbb{N}}$ and each function $f^{j}$, define
scalloped balls $B^{j}_{n}$ and construct a function $\varphi^{j}$
exactly as in Case I above, so that
    $$
    |\varphi^{j}(y)-f^{j}(y)|\leq\varepsilon_{j}
    \hspace{0.3cm}\text{ for all }\hspace{0.3cm}y\in B(y^{j}_{n},
    r^{j}_{n}).
    $$
This function $\varphi^{j}$ is of the form
    $$
    \varphi^{j}(x)=\frac{\sum_{n=1}^{\infty}\alpha^{j}_{n}(x)\varphi^{j}_{n}(x)}
    {\sum_{n=1}^{\infty}\varphi^{j}_{n}(x)}=\lim_{n\to\infty}f^{j}_{n}(x),
    $$
where
    $$
    f^{j}_{n}(x)=\frac{\sum_{k=1}^{n}\alpha^{j}_{k}(x)\varphi^{j}_{k}(x)}
    {\sum_{k=1}^{n}\varphi^{j}_{k}(x)}, \, \textrm{ for all }\,
    x\in \bigcup_{i=1}^{n}B^{j}_{i},
    $$
the domains of the $f^{j}_{n}$ form increasing towers of open sets
whose union is $X$, and, for each $x\in X$ there is some open
neighborhood $V^{j}_{x}$ of $x$ and some $n^{j}_{x}\in\mathbb{N}$
so that $\varphi^{j}(y)=f^{j}_{n}(y)$ for all $y\in V^{j}_{x}$
and all $n\geq n^{j}_{x}$ (see Fact \ref{varphi is locally f_n}).

Now define the mappings $\varphi:X\longrightarrow\mathbb{R}^{m}$
and $f_{n}:\bigcap_{j=1}^{m}\bigcup_{i=1}^{n}B^{j}_{i}
\longrightarrow\mathbb{R}^{m}$ by
    $$
    \varphi(x)=(\varphi^{1}(x),...,\varphi^{m}(x)), \, \text{ and }\,
    f_{n}(x)=(f^{1}_{n}(x), ..., f^{m}_{n}(x)).
    $$
By the choice of the $\varepsilon_j$ and the construction of the
functions $\varphi^{j}$, it is clear that
    $$
    \|\varphi(x)-f(x)\|\leq\varepsilon/2 , \, \textrm{ for all }\,
    x\in X,
    $$
that is, $\varphi$ approximates $f$ as is required.

\begin{fact}\label{the critical set of varphi is locally compact}
If $x\in C_{f_{n}}\cap\big[\bigcap_{j=1}^{m}
\bigcup_{i=1}^{n}B^{j}_{i}\big]$ then $x\in\mathcal{A}[y^{j}_{i}
: 1\leq i\leq n, 1\leq j\leq m]$.
\end{fact}
\begin{proof}
This is a consequence of Fact \ref{the critical set is locally
finite dimensional}. Indeed, according to Remark \ref{the
derivative is a linear combination of x-y_n}, each
$df^{j}_{n}(x)$ is a nontrivial linear combination of the vectors
$(x-y^{j}_k)$ (with $k=1, ..., n$). So, for each $j$ and each
$x\in\bigcup_{i=1}^{n}B^{j}_{i}$ we can assign numbers
$\beta^{j}_{1}(x), ..., \beta^{j}_{n}(x)$ such that at least one
of them does not vanish, and
    $$
    df^{j}_{n}(x)=\sum_{k=1}^{n}\beta^{j}_{k}(x) (x-y^{j}_{k}). \eqno(12)
    $$
Suppose now that $x\in\bigcap_{j=1}^{m}
\bigcup_{i=1}^{n}B^{j}_{i}$ and that the linear map
$df_{n}(x):X\longrightarrow\mathbb{R}^{m}$ is not surjective (that
is, $x$ is a critical point of $f_n$); this means that there are
numbers $\gamma_{1}(x), ..., \gamma_{m}(x)$, not all of them
zero, such that
    $$
    \sum_{j=1}^{m}\gamma_{j}(x)df^{j}_{n}(x)=0. \eqno(13)
    $$
Then, by combining $(12)$ and $(13)$ we get that
    $$
    \sum_{j=1}^{m}\sum_{k=1}^{n}\gamma_{j}(x)\beta^{j}_{k}(x) (x-y^{j}_{k})=0,
    \eqno(14)
    $$
where not all of the numbers $\gamma_{j}(x)\beta^{j}_{k}(x)$
vanish. Since the vectors $y^{j}_{k}$ are all linearly
independent, it follows from $(14)$ that $x$ is in the affine
span of the vectors $y^{j}_{k}$ with $j=1,...,m$; $k=1,...,n$.
\end{proof}

\medskip
As $d\varphi$ is continuous, it is obvious that the set of
critical points $C_{\varphi}$ is closed in $U$. Now we can easily
show that $C_{\varphi}$ is locally compact as well. Indeed, take
$x\in X$. For every $j=1,...,m$ we know that there exists a
neighborhood $V^{j}_{x}$ of $x$ and some $n^{j}_{x}\in\mathbb{N}$
so that $\varphi^{j}(y)=f^{j}_{n}(y)$ for all $y\in V^{j}_{x}$
and every $n\geq n^{j}_{x}$. Fix $n=n_{x}:=\max\{n^{1}_{x}, ...,
n^{m}_{x}\}$, and take $W_{x}$ an open bounded neighborhood of $x$
so that $\overline{W_{x}}\subset V_{x}:=
\bigcap_{j=1}^{m}V^{j}_{x}$. Then we have that
$$\varphi(y)=(\varphi^{1}(y),...,\varphi^{m}(y))=(f^{1}_{n}(y),
..., f^{m}_{n}(y))=f_{n}(y)$$ for all $y\in V_{x}$, and in
particular $V_{x}\subset
\bigcap_{j=1}^{m}\bigcup_{i=1}^{n}B^{j}_{i}.$ Now, according to
Fact \ref{the critical set of varphi is locally compact}, it
follows that $C_{\varphi}\cap V_{x}=C_{f_{n}}\cap V_{x}$ is
contained in an affine subspace of dimension $n m$. In particular
$C_{\varphi}\cap \overline{W_{x}}$ is compact, because it is
closed, bounded, and is contained in a finite-dimensional
subspace. 

\medskip

\begin{rem}\label{remark for e(x)}
{\em Let us say a few words as to the way one has to modify the
above proofs in order to establish Proposition \ref{existence of
varphi} when $\varepsilon$ is a continuous positive function. At
the beginning of the proof of Case I of Proposition \ref{existence
of varphi}, before choosing the $\delta_{x}$, we have to take some
number $\alpha_{x}>0$ so that
$|\varepsilon(y)-\varepsilon(x)|\leq\varepsilon(x)/4$ whenever
$\|y-x\|\leq 2\alpha_{x}$ and then we can find some
$\delta_{x}\leq\alpha_{x}$ so that
$|f(y)-f(x)|\leq\varepsilon(x)/8$ whenever $y\in B(x,
2\delta_{x})$. In particular, after choosing the
$r_{n}=\delta_{x_{n}}$ as in the proof of Case I above, we have
that
    $$
    |f(y)-f(y_{n})|\leq\varepsilon(y_{n})/8, \, \text{ and }\,
    \varepsilon(y_{n})\leq\frac{4}{3}\varepsilon(y) \eqno(15)
    $$
for all $y\in B(y_{n}, r_{n})$. Then we can go on with the proof,
with appropriate modifications, to construct the functions
$\varphi$ and $f_{n}$. Some obvious changes must be made in the
definition of the functions $a_{n}$ and $\alpha_{n}$. Fact
\ref{good approximation} now tells us that
    $$
    |\varphi(y)-f(y_{n})|\leq\varepsilon(y_{n})/4 \eqno(16)
    $$
for all $y\in B(y_{n}, r_{n})$. Then, by combining $(15)$ and
$(16)$ we get that
    $$
    |\varphi(y)-f(y)|\leq |\varphi(y)-f(y_{n})|+|f(y_{n})-f(y)|
    \leq\frac{\varepsilon(y_{n})}{4} + \frac{\varepsilon(y_{n})}{8}=
    \frac{3}{8}\varepsilon(y_{n})\leq\frac{\varepsilon(y)}{2}
    $$
for all $y\in B(y_{n}, r_{n})$ and, since these balls cover $X$,
this proves that $|\varphi(y)-f(y)|\leq\varepsilon(y)/2$ for all
$y\in X$.

In Case II it is enough to define the functions
$\varepsilon_{j}(x)=\varepsilon(x)/\sqrt{4m}$, for $j=1, ..., m$.
The rest of the proof applies just replacing $\varepsilon_{j}$
and $\varepsilon$ with $\varepsilon_{j}(x)$ and $\varepsilon(x)$,
and making some obvious minor modifications as in Case I.}
\end{rem}

\medskip

\begin{center}
{\bf Proof of Theorem \ref{countable union of pairwise disjoint
compact sets of critical points}}
\end{center}

The proof of this result is based on that of Case I of Proposition
\ref{existence of varphi}. We will have to select the numbers
$\lambda_{n}$ with more care, and make sure that the boundaries
of the balls considered have a nice transversality property. An
argument similar to that of Remark \ref{remark for e(x)} shows
that there is no loss of generality in assuming that
$\varepsilon$ is constant.

Suppose that we are at the beginning of the proof of Proposition
\ref{existence of varphi} and we only know that
$X=\bigcup_{n=1}^{\infty}B(y_{n}, s_{n}/2)$, where
$s_{n}=\delta_{x_{n}}$, for some sequence of linearly independent
vectors $(y_{n})$, and
    $$
    |f(y)-f(y_{n})|\leq\varepsilon/4 
    \textrm{ provided } \|y-y_{n}\|\leq \frac{3}{2} s_{n}.
    $$
The following lemma
shows that we can slightly move the radii $s_{n}$ so that, for
any finite selection of centers $y_{n}$, the spheres that are the
boundaries of the balls $B(y_{n}, s_{n})$ have empty intersection
with the affine subspace spanned by those centers.

\begin{lem}\label{controlled intersection of spheres}
We can find a sequence of positive numbers $(r_{n})$ with
$s_{n}\leq r_{n}\leq \frac{3}{2} s_{n}$ so that, if we denote
$S_{n}=\partial B(y_{n}, r_{n})$ then,
\begin{itemize}
\item[{(i)}] for each finite sequence of positive integers
$k_{1}<k_{2}<...<k_{m}$,
    $$
    \mathcal{A}[y_{k_{1}}, ..., y_{k_{m}}]\cap S_{k_{1}}\cap ...\cap
    S_{k_{m}}=\emptyset.
    $$
\item[{(ii)}] for any $n,k\in\mathbb{N}$, $y_{n}\notin S_{k}$.
\end{itemize}
\end{lem}
\begin{proof}
We will define the $r_{n}$ inductively.

For $n=1$ we may take $r_{1}\in [s_{1}, \frac{3}{2}s_{1}]$ so
that $r_{1}$ does not belong to the countable set
$\{\|y_{1}-y_{k}\| : k\in\mathbb{N}\}$; this means that
$y_{k}\notin S_{1}$ for any $k\in\mathbb{N}$. On the other hand,
it is obvious that $\{y_{1}\}\cap S_{1}=\emptyset$.

Assume now that $r_{1}, ..., r_{n}$ have already been chosen in
such a way that the spheres $S_{1}$, ..., $S_{n}$ satisfy $(i)$
and $(ii)$, and let us see how we can find $r_{n+1}$. For any
finite sequence of integers $0<k_{1}< ... <k_{j}\leq n+1$, let us
denote
    $$
    \mathcal{A}_{k_{1},...,k_{j}}=\mathcal{A}[y_{k_{1}}, ...,
    y_{k_{j}}].
    $$
For simplicity, and up to a suitable translation (which obviously
does not affect our problem), we may assume that $y_{n+1}=0$, so
that $\mathcal{A}_{k_{1},...,k_{m},n+1}$ is the $m$-dimensional
vector subspace of $X$ spanned by $y_{k_{1}}$, ..., $y_{k_{m}}$.
Now, for each finite sequence of integers $0<k_{1}<...<k_{m}\leq
n$, consider the map $F_{k_{1}, ...,
k_{m}}:\mathcal{A}_{k_{1},...,k_{m}, n+1}\longrightarrow
\mathbb{R}^{m}$ defined by
    $$
    F_{k_{1}, ..., k_{m}}(x)=\big(\|x-y_{k_{1}}\|^{2}- {r_{k_{1}}}^{2}, ...,
    \|x-y_{k_{m}}\|^{2}- {r_{k_{m}}}^{2}\big).
    $$
Note that
    $$
    D F_{k_{1}, ..., k_{m}}(x)=\big(2(x-y_{k_{1}}), ..., 2(x-y_{k_{m}})\big)
    $$
and therefore $\textrm{rank}\big(D F_{k_{1}, ...,
k_{m}}(x)\big)<m$ if and only if $x\in
\mathcal{A}_{k_{1},...,k_{m}}$. By the induction assumption we
know that
    $$
    S_{k_{1}}\cap ... \cap S_{k_{m}}\cap \mathcal{A}_{k_{1},...,k_{m}}=\emptyset,
    $$
hence it is clear that $\textrm{rank}\big(D F_{k_{1}, ...,
k_{m}}(x)\big)=m$ for all $x\in S_{k_{1}}\cap ... \cap
S_{k_{m}}\cap \mathcal{A}_{k_{1},...,k_{m}, n+1}$. This implies
that
    $$
    M_{k_{1}, ..., k_{m}}:=S_{k_{1}}\cap ... \cap
S_{k_{m}}\cap \mathcal{A}_{k_{1},...,k_{m}, n+1}
    $$
is a compact $m-m=0$-dimensional submanifold of
$\mathcal{A}_{k_{1},...,k_{m}, n+1}$, and in particular
$M_{k_{1}, ..., k_{m}}$ consists of a finite number of points (in
fact two points, but we do not need to know this). Therefore
    $$
    M=\bigcup M_{k_{1}, ..., k_{m}}
    $$
(where the union is taken over all the finite sequences of
integers $0<k_{1}<...<k_{n}\leq n$) is a finite set as well. Now
we have that
    $$I:=\big[s_{n+1}, \frac{3}{2}s_{n+1}\big]\setminus
    \Big(\{\|z\| : z\in M\}\cup\{\|y_{j}\| :
    j\in\mathbb{N}\}\Big)
    $$
is an uncountable subset of the real line, so we can find a number
$r_{n+1}\in I$. With this choice it is clear that
    $$
    S_{k_{1}}\cap ... \cap
S_{k_{m}}\cap S_{n+1}\cap \mathcal{A}_{k_{1},...,k_{m}, n+1}=
M_{k_{1},...,k_{m}}\cap S_{n+1}=\emptyset
    $$
for all finite sequences of integers $0<k_{1}<...<k_{m}<n+1$, and
also $$y_{j}\notin S_{n+1}=\partial B(0, r_{n+1}) \textrm{ for all }
j\in\mathbb{N}.$$ Therefore the spheres $S_{1}$, ..., $S_{n}$,
$S_{n+1}$ satisfy $(i)$ and $(ii)$ as well. By induction the
sequence $(r_{n})$ is thus well defined.
\end{proof}

Now define $B_{n}$, $\varphi_{n}$, $\varphi$, $f_{n}$, as in Case
I of the proof of Proposition \ref{existence of varphi}. All the
properties shown in the proof of \ref{existence of varphi} about
the functions $f_{n}$ and $\varphi$ (in particular Facts
\ref{good approximation} and \ref{the critical set is locally
finite dimensional}) are independent of the way we may choose the
numbers $\lambda_{j}$ in the definitions of $B_j$ and
$\varphi_{j}$. Now we only have to see how we can select those
numbers $\lambda_{j}$ so as to have more control over the set
$C_{\varphi}$ of critical points of $\varphi$ and thus prove the
statement of Theorem \ref{countable union of pairwise disjoint
compact sets of critical points}. We will define the numbers
$\lambda_{n}$ and the open sets $U_{n}$ inductively.

\medskip
\noindent {\bf First step.} Define $\varphi_{1}$ as above and put
$f_{1}(x)=\alpha_{1}(x)$ for all $x\in B_{1}=B(y_{1}, r_{1})$.
Set $\mu_{2}=1/2$, $K_{1}=C_{f_{1}}\cap B_{1}=\{y_{1}\}$, and
$U_{1}=B(y_{1},\mu_{2}r_{1})$.

\medskip
\noindent {\bf Second step.} Fix $\lambda_{2}\in (\mu_{2}, 1)$,
and define $B_{2}$, $\varphi_{2}$, and $f_{2}$ as above.
According to Fact \ref{the critical set is locally finite
dimensional}, we have that
\begin{itemize}
\item[{}] $C_{f_{2}}\cap B_{2}\subset\mathcal{A}[y_{1}, y_{2}]$,
and
\item[{}] $C_{f_{2}}\cap (B_{2}\setminus B_{1})
\subseteq\mathcal{A}[y_{2}]$.
\end{itemize}
We claim that there must exist some $\mu_{3}\in (\lambda_{2}, 1)$
so that $\overline{C_{f_{2}}\cap B_{2}\cap B_{1}}\subset B(y_{1},
\mu_{3}r_{1})$. Otherwise there would exist a sequence $(x_{j})$
in $C_{f_{2}}\cap B_{2}\cap B_{1}$ so that $\|x_{j}-y_{1}\|$ goes
to $r_{1}$ as $j$ goes to $\infty$. Since $C_{f_{2}}\cap
B_{2}\subset\mathcal{A}[y_{1}, y_{2}]$, we may assume, by
compactness, that $x_{j}$ converges to some point
$x_{0}\in\partial B(y_{1}, r_{1})=S_{1}$. If $x_{0}\in B(y_{2},
r_{2})$ then $f_{2}'(x_{0})=0$ (by continuity of $f_{2}'$), and
$x_{0}\neq y_{2}$ (because $y_{2}\notin S_{1}$ by ii) of Lemma
\ref{controlled intersection of spheres}), so
    $$
    f_{2}'(x_{0})=\alpha_{2}'(x_{0})\neq 0,
    $$
a contradiction. Therefore it must be the case that
$x_{0}\in\partial B(y_{2}, r_{2})=S_{2}$. But then
    $$
    x_{0}\in S_{1}\cap S_{2}\cap\mathcal{A}[y_{1}, y_{2}],
    $$
and this contradicts Lemma \ref{controlled intersection of
spheres}.

So let us take $\mu_{3}\in (\lambda_{2}, 1)$ such that
$\overline{C_{f_{2}}\cap B_{2}\cap B_{1}}\subset B(y_{1},
\mu_{3}r_{1})$. Choose also some $\nu_{2}\in
(\mu_{2},\lambda_{2})$. In the case that $y_{2}\in B_{1}$, let us
simply set
\begin{itemize}
\item[{}] $U_{2}=B(y_{2}, r_{2})\cap B(y_{1}, \mu_{3}r_{1})\setminus
\overline{B}(y_{1}, \nu_{2}r_{1})$, and
\item[{}] $K_{2}=\overline{C_{f_{2}}\cap B_{2}\cap B_{1}}\subset
U_{2}$.
\end{itemize}
In the case that $y_{2}\notin B_{1}$, find $\delta_{2}\in (0,
\mu_{3} r_{2})$ so that $B(y_{2}, \delta_{2})\subset
B_{2}\setminus\overline{B_{1}}$, and set
\begin{itemize}
\item[{}] $U_{2}=\big[B(y_{2}, r_{2})\cap B(y_{1}, \mu_{3}r_{1})\setminus
\overline{B}(y_{1}, \nu_{2}r_{1})\big]\cup B(y_{2},\delta_{2})$,
and
\item[{}] $K_{2}=\overline{C_{f_{2}}\cap B_{2}\cap B_{1}}\cup\{y_{2}\}\subset
U_{2}$.
\end{itemize}
Clearly, we have that $C_{f_{2}}\subseteq K_{1}\cup K_{2}$, and
$U_{1}\cap U_{2}=\emptyset$.

\medskip
\noindent {\bf Third step.} Now choose $\lambda_{3}\in (\mu_{3},
1)$ with $\lambda_{3}>1-1/3$, and define $B_{3}$, $\varphi_{3}$,
and $f_{3}$ as above. We have that $f_{3}$ and $f_{2}$ coincide on
$(B_{1}\cup B_{2})\setminus B_{3}$. On $B_{3}$, according to Fact
\ref{the critical set is locally finite dimensional}, we know that
\begin{itemize}
\item[{}] $C_{f_{3}}\cap B_{3}\cap B_{2}\cap B_{1}\subseteq\mathcal{A}[y_{1}, y_{2},
y_{3}]$;
\item[{}] $C_{f_{3}}\cap(B_{3}\cap B_{2}\setminus B_{1})\subseteq\mathcal{A}[y_{2},
y_{3}]$, and $C_{f_{3}}\cap(B_{3}\cap B_{1}\setminus
B_{2})\subseteq\mathcal{A}[y_{1}, y_{3}]$ ;\hfill $(17)$
\item[{}] $C_{f_{3}}\cap(B_{3}\setminus (B_{1}\cup
B_{2}))\subseteq \mathcal{A}[y_{3}]$.
\end{itemize}
Again, there must be some $\mu_{4}\in (\lambda_{3}, 1)$ so that
    $$
    \overline{C_{f_{3}}\cap B_{3}\cap(B_{1}\cup B_{2})}
    \subset B(y_{1}, \mu_{4}r_{1})\cup
    B(y_{2}, \mu_{4}r_{2}).
    $$
Otherwise (bearing in mind the local compactness of
$\mathcal{A}[y_{1}, y_{2}, y_{3}]$), there would exist a sequence
$(x_{j})$ in $C_{f_{3}}\cap B_{3}\cap (B_{1}\cup B_{2})$ so that
$(x_{j})$ converges to some point $x_{0}$ and $(x_{j})$ is not
contained in $B(y_{1}, \mu_{4}r_{1})\cup
    B(y_{2}, \mu_{4}r_{2})$ for any $\mu_{4}<1$. Since a subsequence
of $(x_{j})$ must be contained in one of the sets listed in
$(17)$, we deduce that the limit point $x_{0}$ must belong to one
of the following sets:
\begin{itemize}
\item[{}] $S_{2}\cap S_{1}\cap\mathcal{A}[y_{1}, y_{2},
y_{3}]$;
\item[{}] $S_{2}\cap\mathcal{A}[y_{2},
y_{3}]\setminus B_{1}$;
\item[{}] $S_{1}\cap\mathcal{A}[y_{1},
y_{3}]\setminus B_{2}$,
\end{itemize}
Now we have two cases: either $x_{0}\in B_{3}$, or $x\in\partial
B_{3}$. If $x_{0}\in B_{3}$ then $f_{3}'(x_{0})=0$ (by continuity
of $f_{3}'$), and $x_{0}\neq y_{3}$ (because $y_{3}\notin
S_{1}\cup S_{2}$ by (ii) of Lemma~\ref{controlled intersection of
spheres}), so it follows that
    $$
    f_{3}'(x_{0})=\alpha_{3}'(x_{0})\neq 0,
    $$
a contradiction. On the other hand, if $x_{0}\in\partial B_{3}$
then $x_{0}\in S_{3}$ as well, and now one of the following must
hold:
\begin{itemize}
\item[{}] $x_{0}\in S_{3}\cap S_{2}\cap S_{1}\cap\mathcal{A}[y_{1}, y_{2},
y_{3}]$;
\item[{}] $x_{0}\in S_{3}\cap S_{2}\cap\mathcal{A}[y_{2},
y_{3}]$;
\item[{}] $x_{0}\in S_{3}\cap S_{1}\cap\mathcal{A}[y_{1},
y_{3}]$,
\end{itemize}
but in any case this contradicts Lemma \ref{controlled
intersection of spheres}.

Hence we can take $\mu_{4}\in (\lambda_{3}, 1)$ so that
    $$
    \overline{C_{f_{3}}\cap B_{3}\cap(B_{1}\cup B_{2})}
    \subset B(y_{1}, \mu_{4}r_{1})\cup
    B(y_{2}, \mu_{4}r_{2}).
    $$
Take $\nu_{3}\in (\mu_{3},\lambda_{3})$. Now two possibilities
arise. If $y_{3}\in B_{1}\cup B_{2}$, let us define
    $$
    U_{3}=\Biggl[B(y_{3}, r_{3})\setminus\bigcup_{j=1}^{2}
    \overline{B}(y_{j}, \nu_{3}r_{j})\Biggr]\bigcap
    \Biggl[\bigcup_{j=1}^{2}B(y_{j}, \mu_{4}r_{j})\Biggr],
    $$
and
    $$
    K_{3}=\overline{C_{f_{3}}\cap B_{3}\cap(B_{1}\cup B_{2})}\subset U_{3}.
    $$
If $y_{3}\notin B_{1}\cup B_{2}$, since $y_{3}\notin S_{1}\cup
S_{2}$ we can find $\delta_{3}\in (0, \mu_{4}r_{3})$ so that
$B(y_{3},\delta_{3})\subseteq B_{3}\setminus (B_{1}\cup B_{2})$,
and then we can set
    $$
    U_{3}=\Biggl[\Bigl(B(y_{3}, r_{3})\setminus\bigcup_{j=1}^{2}
    \overline{B}(y_{j}, \nu_{3}r_{j})\Bigr)\bigcap
    \Bigl(\bigcup_{j=1}^{2}B(y_{j}, \mu_{4}r_{j})\Bigr)\Biggr]\bigcup
    B(y_{3}, \delta_{3}),
    $$
and
    $$
    K_{3}=\overline{[C_{f_{3}}\cap B_{3}\cap(B_{1}\cup B_{2})]\cup
    \{y_{3}\}}\subset U_{3}.
    $$
Notice that $U_{3}$ does not meet $U_{1}$ or $U_{2}$, and
$C_{f_{3}}\subseteq K_{1}\cup K_{2}\cup K_{3}$.

\medskip
\noindent {\bf N-th step.} Suppose now that $\mu_{j}$,
$\lambda_{j}$, $\nu_{j}$, $\varphi_{j}$, $B_{j}$, $f_{j}$,
$K_{j}$, $U_{j}$ have already been fixed for $j=1, ..., n$ (and
also $\mu_{n+1}$ has been chosen) in such a manner that $f_{j}$
agrees with $f_{j-1}$ on $(B_{1}\cup ...\cup B_{j-1})\setminus
B_{j}$, and $K_{j}$ and $U_{j}$ are of the form
    $$
    K_{j}=\overline{C_{f_{j}}\cap B_{j}\cap(B_{1}\cup... \cup B_{j-1})}
    \eqno(18)
    $$
and
    $$
    U_{j}=\biggl[B(y_{j}, r_{j})\setminus\Bigl(\bigcup_{i=1}^{j-1}
    \overline{B}(y_{i}, \nu_{j}r_{i})\Bigr)\biggr]\bigcap
    \Bigl[\bigcup_{i=1}^{j-1}B(y_{i}, \mu_{j+1}r_{i})\Bigr] \eqno(19)
    $$
in the case that $y_{j}\in B_{1}\cup ...\cup B_{j-1}$, and are of
this form plus $\{y_{j}\}$ and $B(y_{j},\delta_{j})$ respectively
when $y_{j}\notin B_{1}\cup ...\cup B_{j-1}$; assume additionally
that $U_{j}\cap U_{k}=\emptyset$ whenever $j\neq k$, that
$C_{f_{j}}\subseteq \bigcup_{i=1}^{j}K_{i}$, and that
$\lambda_{j}>1-1/j$. Let us see how we can choose $\lambda_{n+1}$,
$\mu_{n+2}$, $\nu_{n+1}$, $K_{n+1}$ and $U_{n+1}$ so that the
extended bunch keeps the required properties.

Pick any $\lambda_{n+1}\in (\mu_{n+1}, 1)$ so that
$\lambda_{n+1}>1-1/(n+1)$, and define $\varphi_{n+1}$, $B_{n+1}$
and $f_{n+1}$ as above. We know that $f_{n+1}$ agrees with
$f_{n}$ on the set $(B_{1}\cup ... \cup B_{n})\setminus B_{n+1}$.
On $B_{n+1}$, according to Fact \ref{the critical set is locally
finite dimensional}, we have that
    $$
    C_{f_{n+1}}\cap\big(B_{n+1}\setminus\bigcup_{j=1}^{m}B_{k_{j}}\big)
    \subseteq\mathcal{A}\big[\{y_{1}, ..., y_{n+1}\}
    \setminus\{y_{k_{1}}, ..., y_{k_{m}}\}\big]
    $$
for every finite sequence of integers $0<k_{1}<k_{2}< ...
<k_{m}<n+1$.

We claim that there exists some $\mu_{n+2}\in (\lambda_{n+1}, 1)$
so that
    $$
    \overline{C_{f_{n+1}}\cap B_{n+1}\cap (B_{1}\cup ... \cup B_{n})}
    \subseteq \bigcup_{i=1}^{n}B(y_{i},\mu_{n+2}r_{i}).
    $$
Otherwise there would exist a finite (possibly empty!) sequence of
integers $0<k_{1}<k_{2}< ... <k_{m}<n+1$, and a sequence
$(x_{j})_{j=1}^{\infty}$ contained in
    $$
    \Biggl[C_{f_{n+1}}\cap B_{n+1}\cap \Bigl(\bigcap_{j=1}^{\ell}B_{i_{j}}\Bigr)\Biggr]
    \setminus\Bigl(\bigcup_{j=1}^{m}B_{k_{j}}\Bigr)
    \subseteq\mathcal{A}[y_{i_{1}}, ..., y_{i_{\ell}}, y_{n+1}]
    $$
(where $i_{1}, ...,i_{\ell}$ are the positive integers less than
or equal to $n$ that are left when we remove $k_{1}, ...,
k_{m}$), such that $(x_{j})$ converges to some point $x_{0}\in
S_{i_{1}}\cap ...\cap S_{i_{\ell}}$ with
$x_{0}\notin\bigcup_{j=1}^{m}B_{k_{j}}$.

If $x_{0}\in B_{n+1}$ then $f_{n+1}'(x_{0})=0$ (by continuity of
$f_{n+1}'$), and $x_{0}\neq y_{n+1}$, so we easily see that
    $$
    f_{n+1}'(x_{0})=\alpha_{n+1}'(x_{0})\neq 0,
    $$
a contradiction.

If $x_{0}\in \partial B_{n+1}$ then $x_{0}\in S_{n+1}$ as well,
and in this case we have
    $$
    x_{0}\in S_{i_{1}}\cap ...\cap S_{i_{\ell}}\cap S_{n+1}\cap
    \mathcal{A}[y_{i_{1}}, ..., y_{i_{\ell}}, y_{n+1}],
    $$
but this contradicts Lemma \ref{controlled intersection of
spheres}.

Therefore we may take $\mu_{n+2}\in (\lambda_{n+1}, 1)$ so that
    $$
    \overline{C_{f_{n+1}}\cap B_{n+1}\cap (B_{1}\cup... \cup
    B_{n})}\subseteq
    \bigcup_{i=1}^{n}B(y_{i},\mu_{n+2}r_{i}).
    $$
Choose any $\nu_{n+1}\in (\mu_{n+1},\lambda_{n+1})$. As before,
now we face two possibilities. If
$y_{n+1}\in\bigcup_{i=1}^{n}B_{i}$, let us define
    $$
    U_{n+1}=\Biggl[B(y_{n+1}, r_{n+1})\setminus\bigcup_{i=1}^{n}
    \overline{B}(y_{i}, \nu_{n+1}r_{i})\Biggr]\bigcap
    \Biggl[\bigcup_{i=1}^{n}B(y_{i}, \mu_{n+2}r_{i})\Biggr], \,
    \text{ and }
    $$
    $$
    K_{n+1}=\overline{C_{f_{n+1}}\cap B_{n+1}\cap(B_{1}\cup... \cup
    B_{n})}.
    $$
If $y_{n+1}\notin\bigcup_{i=1}^{n}B_{i}$, since $y_{n+1}\notin
S_{i}$ we may find $\delta_{n+1}\in (0, \mu_{n+2}r_{n+1})$ so
that $B(y_{n+1},\delta_{n+1})\subseteq B_{n+1}\setminus
\bigcup_{i=1}^{n}B_{i}$, and then we can add this ball to the
above $U_{n+1}$, and the point $\{y_{n+1}\}$ to that $K_{n+1}$,
in order to obtain sets $U_{n+1}$, $K_{n+1}$ with the required
properties.

By induction, the sequences $(\varphi_{n})$, $(f_{n})$, $(U_{n})$,
$(K_{n})$, $(\lambda_{n})$, $(\mu_{n})$, $(\nu_{n})$ are well
defined and satisfy the above properties. From the construction
it is clear that $U_{n}\cap U_{m}=\emptyset$ whenever $n\neq m$,
and
    $$
    C_{f_{n}}\subseteq\bigcup_{j=1}^{n}K_{j}
    $$
for all $n$. As observed before (see Fact \ref{varphi is locally
f_n}), for each $x\in X$ there exists an open neighborhood
$V_{x}$ of $x$ and some $n_{x}\in\mathbb{N}$ so that
$\varphi(y)=f_{n_{x}}(y)$ for all $y\in V_{x}$. Bearing these
facts in mind, it is immediately checked that
$C_{\varphi}\subseteq\bigcup_{n=1}^{\infty}K_{n}.$ The other
properties in the statement of Theorem \ref{countable union of
pairwise disjoint compact sets of critical points} are immediately
deduced from the above construction.


\medskip

\begin{center}
{\bf Acknowledgements}
\end{center}
\noindent We are indebted to Tadeusz Dobrowolski, who read the
first version of this paper and suggested a simplification of the
original proof by pointing out to us West's paper \cite{West}. We
also wish to thank Pilar Cembranos, Jos\'e Mendoza, Tijani
Pakhrou and Ra\'ul Romero, who helped us to realize that Fact
\ref{the critical set is locally finite dimensional} fails
whenever the norm is not hilbertian. This research was in part
carried out during a stay of the first-named author in the
Mathematics Department of University College London; this author
wish to thank the members of that Department and very especially
David Preiss for their kind hospitality and advice.


\vspace{3mm} \noindent Departamento de An\'alisis Matem\'atico.
Facultad de Ciencias Matem\'aticas. Universidad Complutense.
28040 Madrid, SPAIN\\ \noindent Departamento de An\'alisis
Matem\'atico. Universidad de Sevilla. Sevilla, SPAIN. \noindent
{\em E-mail addresses:} daniel\_azagra@mat.ucm.es, mcb@us.es
\end{document}